\documentclass[times,doublespace]{nlaauth}
\usepackage[dvips,colorlinks,bookmarksopen,bookmarksnumbered,citecolor=red,urlcolor=red]{hyperref}
\usepackage{rotating}
\usepackage{graphicx}
\newtheorem{Thm}{\bf Theorem}[section]

\begin{document}

\runningheads{L Wang, H Liang, FS Bai and Y Huo}{A load balancing strategy for parallel computation of sparse permanents}

\title{A load balancing strategy for parallel computation of sparse permanents}
\author{Lei Wang\affil{1}, Heng Liang\affil{1}'\corrauth'\footnotemark[2],
Fengshan Bai\affil{1} and Yan Huo\affil{2}}

\address{\affilnum{1}Department of Mathematical Sciences,
Tsinghua University, Beijing, 100084, P.R.CHINA.\break
\affilnum{2}China Citic Bank, Block C, Fuhua Mansion, No.8
Chaoyangmen Beidajie, Dongcheng, Beijing, 100027, P.R.CHINA.}

\corraddr{Heng Liang, Department of Mathematical Sciences,
Tsinghua University, Beijing, 100084, P.R.CHINA.}

\cgsn{National Science Foundation of China}{10871115}

\begin{abstract}
The research in parallel machine scheduling in combinatorial
optimization suggests that the desirable parallel efficiency could
be achieved when the jobs are sorted in the non-increasing order of
processing times.  In this paper, we find that the time spending for
computing the permanent of a sparse matrix by hybrid algorithm is
strongly correlated to its permanent value. A strategy  is
introduced to improve a parallel algorithm for sparse permanent.
Methods for approximating permanents, which have been studied
extensively, are used to approximate the permanent values of
sub-matrices to decide the processing order of jobs. This gives an
improved load balancing method. Numerical results show that the
parallel efficiency is improved remarkably for the permanents of
fullerene graphs, which are of great interests in nanoscience.
\end{abstract}


\maketitle
\footnotetext[2]{E-mail: hliang72@gmail.com.}

\section{Introduction}
The permanent of an $n \times n$ matrix $A=[a_{ij}]$ is
defined as
\begin{equation}\label{E1}
per(A)=\sum_{\sigma \in
\Lambda_n}\prod_{i=1}^{n}a_{i\sigma(i)}\hskip 3cm
\end{equation}
where $\Lambda_n$ denotes the set of  all possible permutations of
$\{ 1,2,...,n \} $.

The permanent attracts attentions from mathematics, computer
science, statistical physics and chemical graph theory
\cite{Min78,Tri92}. However, computing the permanent of a matrix is
proved to be a $\#P$-complete problem in counting \cite{Val79},
which is no easier than an $NP$-complete problem in combinatorial optimization.
Even for 3-regular matrices, which are with 3 nonzero entries in each row
and column, evaluating their permanents is still a $\#P$-complete
problem \cite{DL92}.

The best-known algorithm for precise evaluation of the permanent of
general matrix is due to Ryser \cite{Rys63}, and later improved by
Nijenhuis and Wilf \cite{NW78}. It is $O(n2^{n-1})$ in time complexity.
We call the method R-NW algorithm. The R-NW only works for
small matrices.

It is only possible to make the precise calculation faster, if the
special structure properties of matrices can be used intensively.
Several efficient precise algorithms have been proposed by exploring
the structure properties of sparse matrices, such as Kallman's
method \cite{Kal82,DC00}, hybrid algorithm \cite{LB04, LHB06}. Among
them, the hybrid algorithm is the best one for very sparse matrix.

The hybrid algorithm is parallel in nature. A parallelized version
of the algorithm is developed for the permanent computation problem
arising from molecular chemistry \cite{LTB08}. The basic idea of the
parallelized hybrid algorithm is divide and conquer. An $n\times n$
matrix $A$ is divided into a series smaller sub-matrices by using
the hybrid algorithm. When the computational times of the permanents
of sub-matrices are known or estimated appropriately, the load
balancing strategies for the permanent computation could be further
improved with the help of the theory of parallel machine scheduling
in combinatorial optimization. In this paper, we use the statistical
methods to explore the factors which are related to the computational
time of permanent with the hybrid algorithm. An efficient estimation
for computational time of permanent is obtained. Hence the improved
parallel strategy for permanent of sparse graph is proposed.

In the next section, a brief introduction to a hybrid algorithm and
its parallelized version for permanent, which are the best methoths
for very sparse matrix as far as we know, are presented. The load
balancing strategy of parallel algorithm is discussed. In section 3,
the statistical analysis of computational time of permanent is
given. It is shown that the permanent value has strong correlation
to its computational time with the hybrid algorithm. Then an
improved loading balance strategy based on approximate permanent
algorithm is proposed. In section 4, the numerical results are
given. Some discussions are made in section 5.

\section{Parallel algorithms for sparse permanents}
\subsection{A hybrid algorithm for sparse permanent}
Taking the advantage of the sparse structure extensively, a hybrid
method is proposed \cite{LB04, LHB06}. Consider an expansion
\begin{eqnarray}
per \left(\begin{array}{ccccc}a&b&c&d&{\bf x^T}\\{\bf y}_1&{\bf y}_2&{\bf y}_3&{\bf y}_4&{\bf Z}\end{array}\right) \hspace*{-6.0cm} \nonumber\\
 & &=per \left(\begin{array}{cccc}a{\bf y}_2+b{\bf y}_1&{\bf y}_3&{\bf y}_4&{\bf Z}\end{array}\right)
+per \left(\begin{array}{cccc}{\bf y}_1&{\bf y}_2&c{\bf y}_4+d{\bf y}_3&{\bf Z}\end{array}\right) \nonumber\\
 & &\ +per \left(\begin{array}{ccccc}0&0&0&0&{\bf x^T}\\{\bf y}_1&{\bf y}_2&{\bf y}_3&{\bf y}_4&{\bf Z}\end{array}\right) \nonumber\\
 & &=per(A_1)+per(A_2)+per(A_3) \label{E2}
\end{eqnarray}
\noindent where $a$, $b$, $c$ and $d$ are scalars, $x^T$ is an $(n-4)$-dimensional row vector,
${\bf y}_1$, ${\bf y}_2$ ${\bf y}_3$ and ${\bf y}_4$ are both
$(n-1)$-dimensional column vectors, and ${\bf Z}$ is an $(n-1)\times
(n-4)$ matrix. This expansion appears in \cite{FM03}, and is used to
establish an approximate algorithm for permanent. When ``$s < 5$",
where $s$ is the minimal number of nonzero entries in one row or
column of matrix, one
$n \times n$ matrix can be divided into no more than two
$(n-1) \times (n-1)$ matrices, that is, $per(A_3)=0$ in the expansion
(\ref{E2}). Combining the
expansion (\ref{E2}) with R-NW algorithm, a hybrid algorithm is
constructed\cite{LHB06}.

\vskip 3mm
\noindent {\bf Algorithm H$\_$per(Hybrid)}
\begin{description}
  \item[Input:] $A$---an $n\times n$ 0-1 valued matrix.
  \item[Output:] $P=H\_per(A)$.
  \item[Step 1:]  Find the minimal number of nonzero entries
$s$ in one row or column of $A$.
  \item[Step 2:] \begin{description}
                   \item[If] $n>2$ and $s<5$, {\bf then}
                   divide $A$ into $A_1,A_2$ as (\ref{E2}),
                   and $$P=H\_per(A_1)+H\_per(A_2)$$
                   \item[Else] return by $R$-$NW(A)$.
                 \end{description}
\end{description}
It is an efficient algorithm for very sparse matrix, especially
for fullerene-like matrices\cite{LB04, HLB06}.

\subsection{Parallelized version of the hybrid method}
The parallelization of the algorithm is essential for computing
large scale problems. The hybrid algorithm H\_per for permanent is
parallel in nature. First, the $n\times n$ matrix $A$ is divided
into a series of $(n-d)\times (n-d)$ matrices by using the formula
(\ref{E2}) repeatedly. Then the $(n-d)\times (n-d)$ matrices
are computed in parallel. The $d$ is
called the depth of pre-expansion.

Let $A_k^{(w)}$ denote the $w$-th
$(n-k)\times (n-k)$ matrix. Based on the algorithm H\_per, the
following parallel method PH is constructed \cite{LTB08}.

\vskip 3mm
\noindent {\bf Algorithm PH (Parallel H\_per)}
\begin{description}
  \item[Step 1:] Let $n$ be the order of  matrix $A$, num be
the number of  CPU's used, $A_0^{(1)}=A$, $s=1$, set $d$ be the depth of
pre-expansion.
  \item[Step 2:] \begin{description}
                   \item[For]  k=1:d
                   \begin{description}
                     \item[] t=0;
                     \item[] for w=1:s
                     \begin{description}
                       \item[] divide $A_{k-1}^{(w)}$ into $A^{(1)}$
                       and $A^{(2)}$ as (\ref{E2}), $A_k^{(t+1)}=A^{(1)}$, $A_k^{(t+2)}=A^{(2)}$, $t=t+2$;
                     \end{description}
                     \item[] end
                     \item[] s=t;
                   \end{description}
                   \item[End]
                 \end{description}
  \item[Step 3:] Assign $A_d^{(w)}$ in the natural order to the
the currently least load processor and compute the permanent of 
$A_d^{(w)}$ by Algorithm H\_per,
until all $A_d^{(w)}$'s $(1\leq w\leq t)$ has been computed.
  \item[Step 4:] $P=\sum\limits_{w=1}^t H\_per(A_d^{(w)})$.
\end{description}

\subsection{The approximate algorithms for permanents}
Methods for approximating permanents of 0-1 matrices attract a great
deal of studies in the last decade. Markov chain Monte Carlo methods
\cite{JS89,JSV04} absorb great efforts from computer scientists. The
theoretical analysis for those methods are relatively abundant and a
fully-polynomial randomized approximation scheme for the permanent
of arbitrary matrix with non-negative entries has been reported
\cite{JSV04}. But the method is unlikely to be practical in
computations \cite{CRS03}.

A kind of practical approximate methods for permanents is Monte
Carlo method, which reduce permanents to determinants by randomizing
the elements of matrices \cite{CRS03,KKLLL93,GG81}. The idea is
first introduced by Godsil and Gutman \cite{GG81}. It is improved by
Karmarkar et. al. \cite{KKLLL93}, which is one of the most popular
practical approximate algorithms for matrix permanent. Assume
 $$w_0=1,\ \ w_1=-\frac{1}{2}+\frac{\sqrt{3}}{2}i,\ \
w_2=-\frac{1}{2}-\frac{\sqrt{3}}{2}i$$ be the three cube roots of
unity. Let $y$ be a complex number, and $\bar{y}$ denote the
complex conjugate of $y$. The KKLLL method is outlined as follows.

\vskip 3mm
\noindent {\bf Algorithm KKLLL (Karmarkar/Karp/Lipton/Lovasz/Luby)}
\begin{description}
  \item[Input:] $A$---an $n\times n$ 0-1 valued matrix.
  \item[Output:] $X_A$---the estimate for $Per(A)$.
  \item[Step 1:] For all $i,j,\ 1\leq i,j \leq n$,
  \begin{description}
    \item[If] $A_{ij}=0$ then $B_{ij}\leftarrow 0$;
    \item[Elseif] $A_{ij}=1$ then randomly and independently
    choose $B_{ij}\in \{w_0,w_1,w_2\}$ with probability $\frac{1}{3}$.
  \end{description}
  \item[Step 2:] $X_A=det(B)\overline{det(B)}$.
\end{description}
\vskip 3mm

\begin{Thm}
[\cite{KKLLL93}] The KKLLL estimator $X_A$ is unbiased with $E[X_A]=per(A)$.
\end{Thm}

An $(\epsilon,\delta)$-approximation algorithm for $per(A)$ is a
Monte-Carlo algorithm that accepts as input $A$ and two positive
parameters $\epsilon$ and $\delta$. The output of the algorithm
is an estimate $Y$ of $per(A)$, which satisfies
\begin{equation}\label{E3}
P[(1-\epsilon)per(A)\leq Y\leq(1+\epsilon)per(A)]\geq 1-\delta.
\end{equation}

The KKLLL estimator is unbiased and yields an
$(\epsilon,\delta)$-approximation algorithm for estimating $per(A)$
in time $2^{n/2}\frac{1}{\epsilon^2}\log(\frac{1}{\delta}){\bf poly}(n)$
\cite{KKLLL93}. However, for the random 0,1-matrix, Frieze and Jerrum
proved the following result.

\begin{Thm}
[\cite{FJ95}] Let $\omega(n)$ is any function tending to infinity as $n\rightarrow0$.
Then only $O(n\omega(n)\frac{1}{\epsilon^2})$ trials using the KKLLL estimator
suffice to obtain a reliable approximation to the permanent of the random
0,1-matrix within a factor $1\pm\epsilon$ of the correct value.
\end{Thm}

\subsection{The parallel machine scheduling}
The load balancing strategies for the permanent computation can be
further improved with the help of the models of parallel machine
scheduling in combinatorial optimization. Consider the following
machine scheduling model first. Assume that one has a set of $n$
jobs $J_1,\cdots,J_n$, and $m$ identical machines $M_1,\cdots,M_m$.
Each job $J_j$ must be processed without interruption for a time
$p_j>0$ on one of the machines. Each machine can process at most one
job at a time. If all jobs are ready for processing in the very
beginning, it is called {\bf offline} machine scheduling; otherwise
if jobs can only be ready for processing one by one, it is called
{\bf online}.

An algorithm called LS is designed for the online parallel machine
scheduling problems, where jobs are processed in its natural order
of coming. Graham \cite{Gra66} gives the worst-case analysis of the
scheduling heuristics and shows that Algorithm LS has a worst-case
ratio of $2-\frac{1}{m}$, where $m$ is the number of machines
available. If the jobs are sorted in the non-increasing order of
processing times for offline problems, then there is an algorithm
known as LPT. It is proved by Graham that Algorithm LPT has an
improved worst-case ratio of $\frac{4}{3}-\frac{1}{3m}$
\cite{Gra69}.

The scheduling problem in Algorithm PH is essentially online in
which the sub-matrices are sent to the different processors in their
natural order of expansion. In this paper, we will give a
approximate order for the computational times of the sub-matrices.
It is observed and checked that there is a strong correlation
between the permanent value and its computational time. Hence we use
the approximate value of permanent to determine the order of jobs
which are sent to processors so as to improve the parallel
efficiency than the case with natural order.

\section{The improved load balancing strategy}
We find that the time of computing the permanent of a sparse matrix
by hybrid algorithm is strongly correlated to its permanent value.
We also note that the computational time of hybrid algorithm is
dependent on the locations of the nonzero elements of matrix.
Therefore, for any matrix $A=[a_{ij}]_{n\times n}$ construct a 0-1
matrix $B=[b_{ij}]_{n\times n}$ in such a way that  $b_{ij}=1$ if
$a_{ij}\neq0$ and $b_{ij}=0$ if $a_{ij}=0$ for any $1\leq i,j\leq
n$. The computational times of $per(A)$ and $per(B)$ with algorithm
H\_{per} are almost equal. For any matrix $A$, what we discuss is
the relationship between the $per(B)$ and the computational time of
$per(A)$. The following two subsections will present the
statistical analysis for the correlation. Then an improved load
balancing strategy is proposed.

\subsection{Linear regression analysis for computational time of permanent}
The linear regression model is used to investigate which
factors are sensitively response to the computational time of matrix
permanent with algorithm H\_per.

Take the computational time $T$ with algorithm H\_per as the dependent
variable. The following five matrix invariants are considered, which
are chosen empirically and may be related to the computational
time:  the permanent value of the matrix which is denoted as $P$;
the absolute value of the determinant of the matrix which is
denoted as $|D|$; the number of nonzero elements of the matrix
which is denoted as $S$; the variance of sum of nonzero
elements in each row which is denoted as $V_1$; the variance
of sum of nonzero elements in each column which is denoted as $V_2$.

We generate 0-1 matrices randomly with various size and sparsity.
Each group contains 100 matrices. For each 0-1 matrix group, consider the
multiple regression model as follows.
\begin{equation}\label{E4}
\begin{array}{c}
T=a_0+a_1P+a_2|D|+a_3S+a_4V_1+a_5V_2
\end{array}
\end{equation}

The coefficient of determination, often referred to $R^2$, is a frequently used measure
of the fit of the regression line. The definition is, simply,
$$R^2=\frac{\sum_{i=1}^{n}(\hat{y_i}-\bar{y})^2}{\sum_{i=1}^{n}(y_i-\bar{y})^2}$$
corresponding to linear model $y_i=\beta_0+\beta_1X_i+\varepsilon_i (i=1,2,\cdots, n)$,
where $\hat{y}_i=\hat{\beta_0}+\hat{\beta_1}X_i$ and $\bar{y}=\frac{1}{n}\sum_{i=1}^{n}{y_i}$,
$\hat{\beta_0}$ and $\hat{\beta_1}$ being the estimates of $\beta_0$ and $\beta_1$\cite{Mye90}.

The coefficient of determination $R^2$ is no less than 0.8 for every data
group, which is shown in Table~\ref{T1}. In order to find out which
factors significantly correlate to the dependent variable we apply stepwise
regression method to the linear regression models. Using the stepwise
regression, only the permanent value $P$ is significant to the computational
time $T$.

Take the case of n=40 and S=5n as an example. The regression equation
with all five factors is
\begin{equation}\label{E5}
T=39.7636+0.3158\times10^{-7}P-0.8812|D|+0.4175S-
20.0094V_1-12.7279V_2 \end{equation}
with $R^2=0.9190$;

The result of stepwise regression is
\begin{equation}\label{E6}
T=-13.6741+0.3330\times10^{-7}P \end{equation} with $R^2=0.8552$.

$|D|$, $S$, $V_1$ and $V_2$ are very easily computed. However, the
computational time of the permanent can not be predicted by these factors. Only
permanent value itself is strongly correlated to computational time. The result is
reasonable that $per(A)$ is just the number of all the nonzeros
expansion terms in (\ref{E1}), which determines the complexity of the
problem to some extend in nature. Though evaluating $per(A)$ is hard, it is
fortunate that there are many good practical approximate algorithms developed for
permanent.

\subsection{Kendall rank correlation analysis}
For improving the load balancing of the algorithm PH, what is
essentially needed to know is the non-increasing order of the
computational times of all the sub-matrices produced by algorithm PH.

The Kendall rank correlation, also referred to Kendall $\tau$ coefficient,
is a common rank correlation method in the theory of statistical relationship.
This coefficient provides a kind of average measure of the agreement between
two measured quantities. Suppose we have a set of $n$ objects which are being
considered in relation to two properties represented by $x$ and $y$. Numbering
the objects from $1$ to $n$ for the purposes of identification in any order
we please, we may say that they exhibit values $x_1,\ldots,x_n$ according
to $x$ and $y_1,\ldots,y_n$ according to $y$. To any pair of individuals,
say the $i$th and the $j$th$(i<j)$, we will allot an $x$-score, denoted
by $a_{ij}=+1$ if $p_j>p_i$(where $p_i$ is the rank of the $i$th member
according to the $x$-quality) and $a_{ij}=-1$ if $p_j<p_i$, subject only to the
condition that $a_{ij}=-a_{ji}$. Similarly we will allot a $y$-score,
denoted by $b_{ij}$, where $b_{ij}=-b_{ji}$. Denoting $S$ by summing
$a_{ij}b_{ij}$ over all values of $i$ and $j$ from $1$ to $n$, the
Kendall $\tau$ coefficient is defined as\cite{Ken62}:
\begin{equation}\label{E7}
\tau=\frac{\sum_{i<j}{a_{ij}b_{ij}}}{\frac{1}{2}n(n-1)}=\frac{S}{\frac{1}{2}n(n-1)}
\end{equation}

The denominator is the number of pairs of comparison. The Kendall
$\tau$ coefficient have three properties: if the agreement between
the rankings is perfect, i.e. every individual has the same rank in
both, $\tau$ should be $+1$, indicating perfect positive correlation;
if the disagreement is perfect, i.e. one ranking is the inverse of
the other, $\tau$ should be $-1$, indicating perfect negative
correlation; for other arrangements $\tau$ should lie between
these limiting values, and in some acceptable sense increasing
values from $-1$ to $1$ should correspond to increasing agreement
between the ranks.

In practical applications of ranking methods there sometimes arise
cases in which two or more individuals are so similar that no
preference can be expressed between them. The ranking members are
then said to be tied. If there is a tie of $t$ consecutive members
all the scores arising from any pair chosen from them is zero. There
are $\frac{1}{2}t(t-1)$ such pairs. If, therefore, we write
\begin{equation}\label{E8}
T=\frac{1}{2}\sum_t {t(t-1)}
\end{equation}

For ties in one ranking, where $\sum_t$ stands for the summation over
various sets of ties in this ranking, and
\begin{equation}\label{E9}
U=\frac{1}{2}\sum_u {u(u-1)}
\end{equation}

For ties in the other, where $\sum_u$ stands for the summation over
various sets of ties in this ranking, our alternative form of the
coefficient $\tau$ for tied ranks may be written
\begin{equation}\label{E10}
\tau=\frac{S}{\sqrt{\frac{1}{2}n(n-1)-T}\sqrt{\frac{1}{2}n(n-1)-U}}
\end{equation}

We can use the distribution of the Kendall $\tau$ coefficient in testing
the significance of $\tau$ under the null hypothesis that the two qualities
are independent. In the null hypothesis case the exact distribution of
$\tau$ can be calculated exactly for small samples, and as $n$ increases
it has been proved that the distribution tends to normality, with
$E(\tau)=0$ and $var(\tau)=\frac{2(2n+5)}{9n(n-1)}$\cite{Ken62}.

We use Kendall rank correlation to measure the association between
the ranking of the computational times by algorithm
H\_per and the rankings by the permanents and approximate permanents.

For a set of 0-1 matrices $A_1,\cdots,A_m$, let the $T_i$ denote the
computational time of the algorithm H\_per, $P_i$ denote the exact
value of $per(A_i)$, $AP_i$ denote the approximate value of
$per(A_i)$ by KKLLL algorithm, $i=1,\cdots,m$.

The similarities between rank of $\{T_i\}$ and rank of $\{P_i\}$,
between rank of $\{T_i\}$ and rank of $\{AP_i\}$ are considered
respectively.

{\noindent \bf Case study 1: Kendall rank correlation analysis for
random 0-1 matrix:} Tested by the Kendall rank correlation coefficient, the
similarities between the ranking of $\{T_i\}$ and that of
$\{P_i\}$, between the ranking of $\{T_i\}$ and that of $\{AP_i\}$ are
significant for all the matrix groups used in subsection 3.1. Take
$(n,S)=(60,4n)$ as an example, the Table~\ref{T2} shows the Kendall
$\tau$ coefficients and $p$-values. In the Kendall rank corelation analysis,
the $p$-value is used for testing
the significance of $\tau$ under the null hypothesis that the two qualities
are independent against the alternative that the two qualities are
dependent. If the $p$-value is small, say less than $0.05$, then the
two qualities are significantly dependent.

{\noindent \bf Case stduy 2: Kendall rank correlation analysis for
3-regular matrix:} The second example comes from chemical graph
theory. Consider the adjacent matrix of a fullerene with 100 atoms,
which is a 3-regular $100\times 100$ matrix. It is divided into
159 $80\times80$ matrices by using the formula (\ref{E2})
repeatedly. The Table~\ref{T3}
illustrates the result of Kendall $\tau$ test.

{\noindent \bf Case study 3: Kendall rank correlation analysis for
4-regular matrix:} The third example is computing $per(I+A)$,
where $A$ is the adjacent matrix of
buckministerfullerene $C_{60}$, $I$ is identity matrix.
The matrix $I+A$ is 4-regular. The $60\times60$ matrix $I+A$ is divided
into 123 $52\times52$ matrices by using the formula (\ref{E2})
repeatedly. The Table~\ref{T4}
illustrates the result of Kendall $\tau$ test.

The p-values in the three cases are all extremely small. Hence the both
rankings of $\{P_i\}$ and $\{AP_i\}$ and that of  $\{T_i\}$
are dependent significantly. The results show that the
computational time with Algorithm H\_per has a strong rank
correlation with the permanent value.

\subsection{The improved load balancing strategy}
The parallel algorithm is improved by taking advantage
of the ordering of estimated permanents as load balancing
strategy. For the algorithm PH, the step 3 is changed as follows.

\vskip 3mm
\begin{description}
  \item[Step 3:] approximate the permanents of $A_d^{(w)} (w=1,2,\cdots,t)$ by KKLLL algorithm, sort the matrices $A_d^{(w)} (w=1,2,\cdots,t)$ as the non-increasing order of the approximate permanents, then assign $A_d^{(w)}$ in the sorted order to the currently least loaded processor and compute the permanent of $A_d^{(w)}$ by Algorithm H\_per, until all $A_d^{(w)}$'s $(1\leq w\leq t)$ has been computed.
\end{description}

\section{Numerical results}
We use the approximate permanent to give the approximate order with
which the sub-matrices divided by algorithm PH are sent to
processors. In this section, the performance of the improved load
balancing strategy is tested by the numerical examples, which are
arising from molecular chemistry application and choosing from sparse
matrix collection. All numerical
experiments in this paper are carried on a 32-bit Intel Pentium III
(1266 MH) with 32 processors, and the programming language is
Fortran 90.

\subsection{The numerical result for the permanent of $C_{100}$ }
The example of fullerene $C_{100}$ used in subsection 3.2 is
considered again. The adjacent matrix $A$ of $C_{100}$ is divided
into 159 sub-matrices. These sub-matrices are sent to the processors
with the three order strategies. One is their natural order of
expansion, which is the strategy of PH algorithm. The second is
the non-increasing order of the exact computational times of
the sub-matrices, which is the ideal strategy according to parallel
machine scheduling. The third is the non-increasing order of
the approximate permanent values of the sub-matrices by KKLLL
algorithm, which is the improved strategy proposed in subsection 3.3.
Approximate permanents play a role of the preconditioning,
which consumes only a little time compared with the computational
time of permanent.

\begin{figure}
\centering
  \includegraphics[width=10cm,height=8cm]{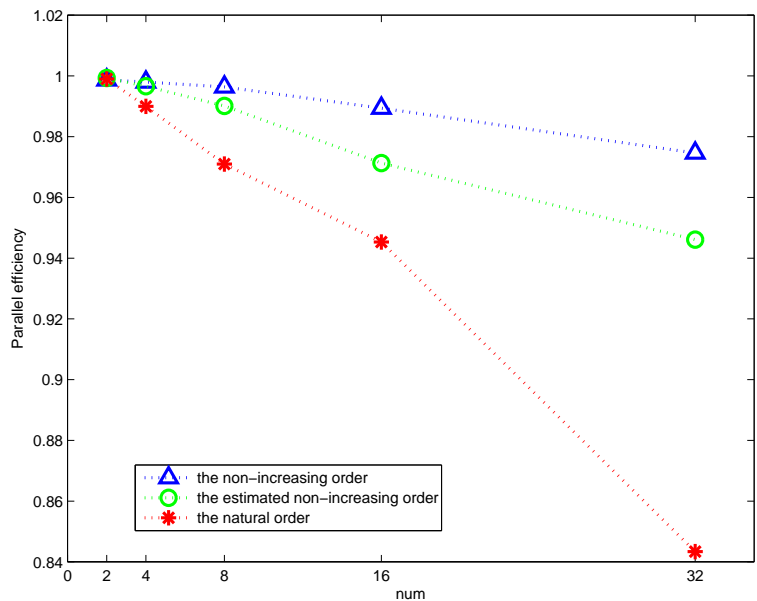}\\
  \caption{The comparison of parallel efficiency for $C_{100}$}
  \label{F1}
\end{figure}

The numerical results with the three order strategies are
shown in the Tables~\ref{T5}-\ref{T7}. The parallel efficiencies
of three order strategies are compared in Figure \ref{F1}.
The results of improved strategy
in Table~\ref{T7} are all better than that of the natural order in
Table~\ref{T5}. Moreover the parallel efficiency in Table~\ref{T7} is
almost the same with that in Table~\ref{T6} except the case of
32 CPU's. But the efficiency of 32 CPU's has reached 94.61\%,
which has been good enough in parallel computation.

\subsection{The numerical result for permanental polynomial of $C_{60}$ }
The permanental polynomial of a graph $G$ is of interest in chemical
graph theory \cite{Tri92}. It is defined as
\begin{equation}\label{E11}
 P(G,x)=per(xI-A),
\end{equation}
where $A$ is the adjacency matrix of the graph $G$ with $n$
vertices, and $I$ is the identity matrix of order $n$. The
permanental polynomial can be obtained by a series of
computations of the permanents formed $per(xI-A)$, where
$x$ is one of the $(n+1)$-th roots of unity in complex plane
\cite{HLB06}. The permanental polynomial of Buckminsterfullerene
$C_{60}$ is first computed by parallel algorithm PH \cite{LTB08}.

\begin{figure}
\centering
  \includegraphics[width=10cm,height=8cm]{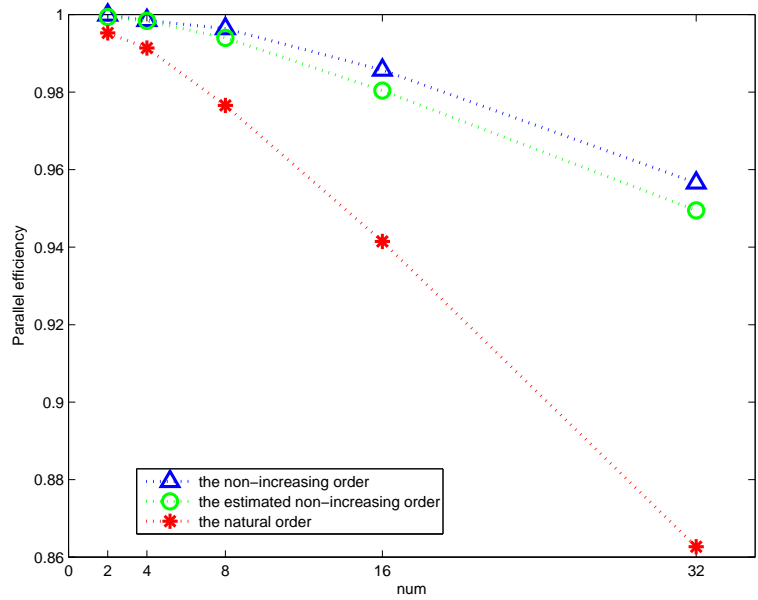}\\
  \caption{The comparison of parallel efficiency for $C_{60}$}
  \label{F2}
\end{figure}

The results of natural order, non-increasing order and estimated
non-increasing order are shown in Tables~\ref{T8}-\ref{T10} respectively.
The parallel efficiencies of three order strategies are compared
in Figure \ref{F2}.
The parallel efficiency under exact non-increasing order is 95.66\%
for 32 CPU's, while that of improved strategy in estimated
non-increasing order is 94.95\%.

\subsection{The numerical result for the permanent of sparse matrix}
We choose a sparse matrix from the University of Florida Sparse Matrix
Collection\cite{CISE} as our third example, which is a large and actively
growing set of sparse matrices that arise in real applications.
This symmetric $66\times 66$ matrix $B$ named $dwt\_66$ has $320$ nonzero elements, and
comes from symmetric connection table from DTNSRDC(David W. Taylor Naval
Ship Research and Development Center), WASHINGTON\cite{CISE}.
We divide it into $140$ $53\times 53$ sub-matrices. These sub-matrices
are sent to the processors
with the three order strategies: the natural order of
expansion, the non-increasing order of the exact computational times of
the sub-matrices, the non-increasing order of
the approximate permanent values of the sub-matrices by KKLLL
algorithm.

\begin{figure}
\centering
  \includegraphics[width=10cm,height=8cm]{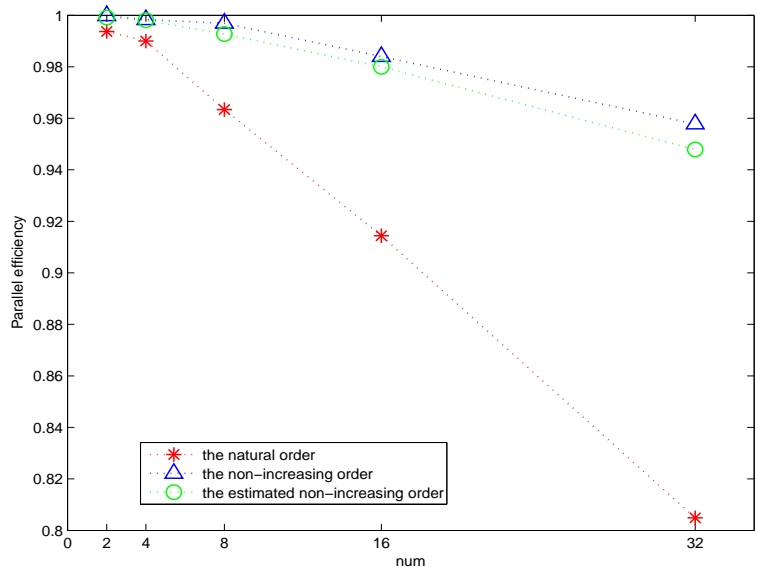}\\
  \caption{The comparison of parallel efficiency for $B$}
  \label{F3}
\end{figure}

The numerical results with the three order strategies are
shown in the Tables~\ref{T15}-\ref{T17}. The parallel efficiencies
of three order strategies are compared in Figure \ref{F3}.
The results of improved strategy
in Table~\ref{T17} are all better than that of the natural order in
Table~\ref{T15}. Moreover the parallel efficiency in Table~\ref{T17} is
almost the same with that in Table~\ref{T16} except the case of
32 CPU's. But the efficiency of 32 CPU's has reached 94.79\%,
which has been good enough in parallel computation comparing with
the efficiency of the natural order 80.49\%.

We have showed the accelerated ratio and parallel efficiency for the three numerical
experiments of parallel algorithm $PH$ using three different load balancing strategies.
The results roughly increase $10\%-15\%$ parallel efficiency using the improved strategy
in contrast to the natural order strategy. Also the parallel efficiency using the
improved strategy is very close to the one using the increasing order strategy,
which is the optimal load balancing strategy for parallel computing.
The permanent values of groups of sub-matrices in our experiments are
given in Appendices.

\section{Conclusion and discussion}
The matrix permanent has critical applications in combinatorial
counting, statistical physics and molecular chemistry. For large
scale matrices,  parallel methods are developing quickly in recent
years. From the results  of parallel machine scheduling  in
combinatorial optimization, one knows that the desirable parallel
efficiency will be achieved when the jobs are sorted in the
non-increasing order of their processing times. Hence it is desired
to know the processing times of the jobs for achieving good parallel
efficiency. In this paper, we find that there are strong correlation
between the permanent value and its computational time. Therefore
the approximate algorithms for permanent are used to estimate the
computational times of sub-matrices, which are the jobs in the
permanent parallel algorithm. The numerical experiments on
fullerene-type graphs, which are of great interest in fullerene
chemistry, show that the parallel efficiency is improved remarkably
by our  load balancing strategy.

The approximate method for matrix permanent used in the paper can
also be regarded as a preconditioner for the parallel hybrid
algorithm\cite{LTB08}. Preconditioning is so successful and valuable
in numerical linear algebra. Following the similar idea, it is
meaningful to establish the basic concepts and a general framework
of the precondition methods for permanent computation, by deeply
investigating the mechanism of the existing successful algorithms.
It is our future work to develop the preconditions such that  the
efficiency of the algorithms for permanents can be highly improved,
and the realistic scientific computation problems can be solved.

\section*{Appendices}
The adjacent matrix of $C_{100}$ (written in MATLAB) is given below, and the permanent values of
159 sub-matrices for matrix $A$ of $C_{100}$ is listed in Table~\ref{T12}.

\begin{tabular}{ccccccccccccccccc}
  D=[ & 2 & 4 & 8; & 7 & 22 & 23; & 23 & 42 & 43; & 43 & 62 & 63; & 63 & 82 & 85; & $\cdots$ \\
   & 1 & 3 & 11; & 19 & 21 & 39; & 39 & 41 & 59; & 59 & 61 & 79; & 79 & 81 & 94; & $\cdots$ \\
   & 2 & 5 & 14; & 21 & 26 & 41; & 41 & 46 & 61; & 61 & 66 & 81; & 66 & 84 & 85; & $\cdots$ \\
   & 1 & 5 & 6; & 9 & 25 & 26; & 26 & 45 & 46; & 46 & 65 & 66; & 67 & 83 & 88; & $\cdots$ \\
   & 3 & 4 & 17; & 10 & 24 & 27; & 27 & 44 & 47; & 47 & 64 & 67; & 81 & 83 & 96; & $\cdots$ \\
   & 4 & 7 & 19; & 23 & 24 & 44; & 43 & 44 & 64; & 63 & 64 & 83; & 70 & 87 & 88; & $\cdots$ \\
   & 6 & 9 & 21; & 25 & 30 & 45; & 45 & 50 & 65; & 65 & 70 & 84; & 71 & 86 & 91; & $\cdots$ \\
   & 1 & 9 & 10; & 12 & 29 & 30; & 30 & 49 & 50; & 50 & 69 & 70; & 84 & 86 & 97; & $\cdots$ \\
   & 7 & 8 & 24; & 13 & 28 & 31; & 31 & 48 & 51; & 51 & 68 & 71; & 74 & 90 & 91; & $\cdots$ \\
   & 8 & 12 & 25; & 27 & 28 & 48; & 47 & 48 & 68; & 67 & 68 & 86; & 75 & 89 & 95; & $\cdots$ \\
   & 2 & 12 & 13; & 29 & 34 & 49; & 49 & 54 & 69; & 69 & 74 & 87; & 87 & 89 & 98; & $\cdots$ \\
   & 10 & 11 & 28; & 15 & 33 & 34; & 34 & 53 & 54; & 54 & 73 & 74 & 78 & 93 & 95; & $\cdots$ \\
   & 11 & 15 & 29; & 16 & 32 & 35; & 35 & 52 & 55; & 55 & 72 & 75; & 80 & 92 & 94; & $\cdots$ \\
   & 3 & 15 & 16; & 31 & 32 & 52; & 51 & 52 & 72; & 71 & 72 & 89; & 82 & 93 & 99; & $\cdots$ \\
   & 13 & 14 & 32; & 33 & 38 & 53; & 53 & 58 & 73; & 73 & 78 & 90; & 90 & 92 & 100; & $\cdots$ \\
   & 14 & 18 & 33; & 18 & 37 & 38; & 38 & 57 & 58; & 58 & 77 & 78; & 85 & 97 & 99; & $\cdots$ \\
   & 5 & 18 & 20; & 20 & 36 & 40; & 40 & 56 & 60; & 60 & 76 & 80; & 88 & 96 & 98; & $\cdots$ \\
   & 16 & 17 & 36; & 35 & 36 & 56; & 55 & 56 & 76; & 75 & 76 & 92; & 91 & 97 & 100; & $\cdots$ \\
   & 6 & 20 & 22; & 22 & 40 & 42; & 42 & 60 & 62; & 62 & 80 & 82; & 94 & 96 & 100; & $\cdots$ \\
   & 17 & 19 & 37; & 37 & 39 & 57; & 57 & 59 & 77; & 77 & 79 & 93; & 95 & 98 & 99; & ] \\
\end{tabular}

$A=zeros(100);\ for\ k=1:100\ A(k,D(k,:))=1;\ end$

\vskip 3mm
\noindent The adjacent matrix of $C_{60}$ (written in MATLAB) is given below, and the permanent values of
123 sub-matrices for matrix $B=xI-A$ of $C_{60}$ is listed in Table~\ref{T14}.

\begin{tabular}{ccccccccccccccccc}
  D=[ & 2 & 3 & 12; & 6 & 12 & 16; & 26 & 27 & 36; & 30 & 36 & 40; & 50 & 52 & 55; & $\cdots$ \\
   & 1 & 4 & 6; & 10 & 15 & 21; & 25 & 28 & 30; & 5 & 34 & 39; & 15 & 49 & 51; & $\cdots$ \\
   & 1 & 11 & 35; & 14 & 16 & 50; & 19 & 25 & 35; & 38 & 40 & 58; & 24 & 45 & 50; & $\cdots$ \\
   & 2 & 5 & 33; & 13 & 15 & 55; & 17 & 26 & 29; & 37 & 39 & 43; & 48 & 49 & 56; & $\cdots$ \\
   & 4 & 8 & 38; & 18 & 19 & 28; & 22 & 28 & 32; & 42 & 44 & 47; & 54 & 56 & 59; & $\cdots$ \\
   & 2 & 7 & 13; & 17 & 20 & 22; & 26 & 31 & 37; & 31 & 41 & 43; & 7 & 53 & 55; & $\cdots$ \\
   & 6 & 8 & 54; & 11 & 17 & 27; & 30 & 32 & 42; & 40 & 42 & 57; & 16 & 49 & 54; & $\cdots$ \\
   & 5 & 7 & 59; & 9 & 18 & 21; & 29 & 31 & 47; & 41 & 48 & 60; & 52 & 53 & 60; & $\cdots$ \\
   & 10 & 11 & 20; & 14 & 20 & 24; & 4 & 34 & 35; & 46 & 48 & 51; & 43 & 58 & 60; & $\cdots$ \\
   & 9 & 12 & 14; & 18 & 23 & 29; & 33 & 36 & 38; & 23 & 45 & 47; & 39 & 57 & 59; & $\cdots$ \\
   & 3 & 9 & 19; & 22 & 24 & 46; & 3 & 27 & 33; & 32 & 41 & 46; & 8 & 53 & 58; & $\cdots$ \\
   & 1 & 10 & 13; & 21 & 23 & 51; & 25 & 34 & 37; & 44 & 45 & 52; & 44 & 56 & 57; & ] \\
\end{tabular}

$A=zeros(60);\ for\ k=1:60\ A(k,D(k,:))=1;\ end$

$x\footnotemark[3]=0.9947+0.1028i;\ B=x*eye(60)-A;$
\footnotetext[3]{We select one of the $61$-th roots of unity in complex plane as $x$.}

\newpage

\begin{table}
\caption{$R^2$ of model for various $n$ and $S$.\label{T1}}
\centering
\tabsize
\begin{tabular}{c|cccc}
\toprule $n$ & $S=4n$ & $S=5n$ & $S=6n$ & $S=7n$\\
\midrule
20 & 0.8822 & 0.8449 & 0.8080 & 0.8296 \\
25 & 0.9110 & 0.8003 & 0.8098 & 0.8713 \\
30 & 0.8916 & 0.8015 & 0.8310 & 0.8620 \\
35 & 0.9451 & 0.8786 & 0.8537 \\
40 & 0.9929 & 0.9190 & 0.8681 \\
45 & 0.8950 & 0.8824 \\
50 & 0.8200 & 0.8902 \\
55 & 0.9978 \\
60 & 0.9639 \\
\bottomrule
\end{tabular}
\end{table}

\begin{table}
\caption{Kendall $\tau$ rank correlation for random matrix with $n=60, S=4n$.\label{T2}}
\centering
\tabsize
\begin{tabular}{c|cc}
\toprule $ $ & coefficient & p-value\\
\midrule
$\{T_{i}\}$ and $\{P_{i}\}$ & 0.4457 & $9.6723\times10^{-17}$ \\
$\{T_{i}\}$ and $\{AP_{i}\}$ & 0.4922 & $3.4198\times10^{-19}$ \\
\bottomrule
\end{tabular}
\end{table}

\begin{table}
\caption{Kendall $\tau$ rank correlation for $C_{100}$.\label{T3}}
\centering
\tabsize
\begin{tabular}{c|cc}
\toprule $ $ & coefficient & p-value\\
\midrule
$\{T_{i}\}$ and $\{P_{i}\}$ & 0.5334 & $2.0294\times10^{-23}$ \\
$\{T_{i}\}$ and $\{AP_{i}\}$ & 0.4919 & $3.6773\times10^{-20}$ \\
\bottomrule
\end{tabular}
\end{table}

\begin{table}
\caption{Kendall $\tau$ rank correlation for $I+A$.\label{T4}}
\centering
\tabsize
\begin{tabular}{c|cc}
\toprule $ $ & coefficient & p-value\\
\midrule
$\{T_{i}\}$ and $\{P_{i}\}$ & 0.4475 & $4.6284\times10^{-17}$ \\
$\{T_{i}\}$ and $\{AP_{i}\}$ & 0.5106 & $1.6604\times10^{-21}$ \\
\bottomrule
\end{tabular}
\end{table}

\begin{table}
\caption{Results of \textit{the natural order} for $per(A)$.\label{T5}}
\centering
\tabsize
\begin{tabular}{c|ccc}
\toprule num & Time(sec) & Accelerated ratio & Parallel efficiency\\
\midrule
1 & 118413.21\footnotemark[1] & -- & -- \\
2 & 59265.87 & 1.99 & 0.9990 \\
4 & 29902.32 & 3.96 & 0.9900 \\
8 & 15243.71 & 7.77 & 0.9710 \\
16 & 7829.07 & 15.13 & 0.9453 \\
32 & 4387.49 & 26.99 & 0.8434 \\
\bottomrule
\end{tabular}
\end{table}

\footnotetext[1]{It takes about $35$ hours to compute $159$ sub-matrices of $A$
of $C_{100}$, and the longest one takes about $40$ minutes while the shortest one takes
about $4$ minutes. While the time to estimate them by the $KKLLL$ algorithm(the trials
is about $n^2$) is extremely short with only $15$ seconds.}

\begin{table}
\caption{Results of \textit{the non-increasing order} for $per(A)$.\label{T6}}
\centering
\tabsize
\begin{tabular}{c|ccc}
\toprule num & Time(sec) & Accelerated ratio & Parallel efficiency\\
\midrule
1 & 118413.21 & -- & -- \\
2 & 59283.67 & 1.99 & 0.9987 \\
4 & 29662.62 & 3.99 & 0.9980 \\
8 & 14855.12 & 7.97 & 0.9964 \\
16 & 7480.11 & 15.83 & 0.9894 \\
32 & 3796.85 & 31.19 & 0.9746 \\
\bottomrule
\end{tabular}
\end{table}

\begin{table}
\caption{Results of \textit{the estimated non-increasing order for $per(A)$}.\label{T7}}
\centering
\tabsize
\begin{tabular}{c|ccc}
\toprule num & Time(sec) & Accelerated ratio & Parallel efficiency\\
\midrule
1 & 118413.21 & -- & -- \\
2 & 59248.07 & 1.99 & 0.9993 \\
4 & 29704.29 & 3.99 & 0.9966 \\
8 & 14949.65 & 7.92 & 0.9901 \\
16 & 7619.50 & 15.54 & 0.9713 \\
32 & 3911.22 & 29.32 & 0.9461 \\
\bottomrule
\end{tabular}
\end{table}

\begin{table}
\caption{Results of \textit{the natural order} for $per(xI-A)$.\label{T8}}
\centering
\tabsize
\begin{tabular}{c|ccc}
\toprule num & Time(sec) & Accelerated ratio & Parallel efficiency\\
\midrule
1 & 125528.59\footnotemark[2] & -- & -- \\
2 & 63060.68 & 1.99 & 0.9953 \\
4 & 31654.37 & 3.97 & 0.9914 \\
8 & 16067.04 & 7.81 & 0.9766 \\
16 & 8333.01 & 15.06 & 0.9415 \\
32 & 4547.08 & 27.61 & 0.8627 \\
\bottomrule
\end{tabular}
\end{table}

\footnotetext[2]{For the $123$ sub-matrices of $xI-A$ of $C_{60}$,
the total time is about $33$ hours as $60$ minutes being the longest
time and $3$ minutes being the shortest one. Compared to the time to
compute this group of sub-matrices, the time to estimate them by the
$KKLLL$ algorithm(the trials is about $n^2$) is extremely short with
only $10$ seconds.}

\begin{table}
\caption{Results of \textit{the non-increasing order} for $per(xI-A)$.\label{T9}}
\centering
\tabsize
\begin{tabular}{c|ccc}
\toprule num & Time(sec) & Accelerated ratio & Parallel efficiency\\
\midrule
1 & 125528.59 & -- & -- \\
2 & 62771.94 & 1.99 & 0.9999 \\
4 & 31430.37 & 3.98 & 0.9985 \\
8 & 15748.20 & 7.93 & 0.9964 \\
16 & 7959.24 & 15.70 & 0.9857 \\
32 & 4100.71 & 30.46 & 0.9566 \\
\bottomrule
\end{tabular}
\end{table}

\begin{table}
\caption{Results of \textit{the estimated non-increasing order for} $per(xI-A)$. \label{T10}}
\centering
\tabsize
\begin{tabular}{c|ccc}
\toprule num & Time(sec) & Accelerated ratio & Parallel efficiency\\
\midrule
1 & 125528.59 & -- & -- \\
2 & 62801.97 & 1.99 & 0.9994 \\
4 & 31432.43 & 3.99 & 0.9984 \\
8 & 15785.78 & 7.95 & 0.9940 \\
16 & 8002.38 & 15.69 & 0.9804 \\
32 & 4131.40 & 30.38 & 0.9495 \\
\bottomrule
\end{tabular}
\end{table}

\begin{table}
\caption{Results of \textit{the natural order} for $per(B)$.\label{T15}}
\centering
\tabsize
\begin{tabular}{c|ccc}
\toprule num & Time(sec) & Accelerated ratio & Parallel efficiency\\
\midrule
1 & 2170843.91\footnotemark[3] & -- & -- \\
2 & 1090876.33 & 1.99 & 0.9937 \\
4 & 548192.90 & 3.96 & 0.9900 \\
8 & 281562.11 & 7.71 & 0.9634 \\
16 & 148383.04 & 14.63 & 0.9144 \\
32 & 84271.89 & 25.76 & 0.8049 \\
\bottomrule
\end{tabular}
\end{table}

\footnotetext[3]{For the $140$ sub-matrices of $B$,
the total time is about $25$ days as $17$ hours being the longest
time and $54$ minutes being the shortest one. Compared to the time to
compute this group of sub-matrices, the time to estimate them by the
$KKLLL$ algorithm(the trials is about $n^2$) is extremely short with
only $20$ minutes.}

\begin{table}
\caption{Results of \textit{the non-increasing order} for $per(B)$.\label{T16}}
\centering
\tabsize
\begin{tabular}{c|ccc}
\toprule num & Time(sec) & Accelerated ratio & Parallel efficiency\\
\midrule
1 & 2170843.91 & -- & -- \\
2 & 1090763.83 & 1.99 & 0.9998 \\
4 & 544071.15 & 3.99 & 0.9983 \\
8 & 272035.57 & 7.98 & 0.9969 \\
16 & 137918.92 & 15.74 & 0.9840 \\
32 & 70826.88 & 30.65 & 0.9578 \\
\bottomrule
\end{tabular}
\end{table}

\begin{table}
\caption{Results of \textit{the estimated non-increasing order for} $per(B)$. \label{T17}}
\centering
\tabsize
\begin{tabular}{c|ccc}
\toprule num & Time(sec) & Accelerated ratio & Parallel efficiency\\
\midrule
1 & 2170843.91 & -- & -- \\
2 & 1090796.83 & 1.99 & 0.9991 \\
4 & 544082.47 & 3.99 & 0.9981 \\
8 & 273406.03 & 7.94 & 0.9927 \\
16 & 138446.67 & 15.68 & 0.9800 \\
32 & 71574.14 & 30.33 & 0.9479 \\
\bottomrule
\end{tabular}
\end{table}

\begin{table}
\centering
\caption{The permanent values of 159 sub-matrices for $C_{100}$. \label{T12}}
\begin{tabular}{cccccccc}
\hline
$No.$ & $per$ & $No.$ & $per$ & $No.$ & $per$ & $No.$ & $per$\\
\hline
1 & 1211353365376 & 41 & 922687029451 & 81 & 410635523677 & 121 & 384774200374\\
2 & 839305965468 & 42 & 427234459851 & 82 & 791821507266 & 122 & 182186664250\\
3 & 965650812160 & 43 & 1294043856556 & 83 & 434949320590 & 123 & 1809089802850\\
4 & 1127638304250 & 44 & 523754341006 & 84 & 852715937662 & 124 & 1013388601052\\
5 & 613075958258 & 45 & 1241484191464 & 85 & 646478774752 & 125 & 703618750146\\
6 & 837748730198 & 46 & 693179603258 & 86 & 633469462226 & 126 & 1239431737028\\
7 & 1127638304250 & 47 & 989214913602 & 87 & 1060604731392 & 127 & 634095267912\\
8 & 555687615100 & 48 & 404136166254 & 88 & 1035132631948 & 128 & 823519791998\\
9 & 895137073356 & 49 & 728857315908 & 89 & 891202457230 & 129 & 1181911459556\\
10 & 802594514720 & 50 & 579851490378 & 90 & 637500159666 & 130 & 658941312496\\
11 & 543460049516 & 51 & 1330651227228 & 91 & 1562498894180 & 131 & 1245530450108\\
12 & 1561971554400 & 52 & 738531936266 & 92 & 670707808150 & 132 & 874154628868\\
13 & 1191963356590 & 53 & 787517656343 & 93 & 797583809616 & 133 & 980918711071\\
14 & 1296946798176 & 54 & 1988491607655 & 94 & 1312280474496 & 134 & 402430987203\\
15 & 885422161590 & 55 & 913255188514 & 95 & 341597043970 & 135 & 749314476366\\
16 & 553665911338 & 56 & 501759364222 & 96 & 751516266132 & 136 & 1557977902746\\
17 & 753899681772 & 57 & 516942933654 & 97 & 1239431737028 & 137 & 282509587768\\
18 & 1396998591510 & 58 & 1597921910139 & 98 & 634095267912 & 138 & 1060604731392\\
19 & 482000132650 & 59 & 683683805805 & 99 & 823519791998 & 139 & 723347758330\\
20 & 1464052439574 & 60 & 1605123039194 & 100 & 676166974426 & 140 & 998475895154\\
21 & 1222801926728 & 61 & 1043139793514 & 101 & 1181911459556 & 141 & 379478870110\\
22 & 588614405964 & 62 & 1213481335486 & 102 & 677714670328 & 142 & 532132377232\\
23 & 731981344504 & 63 & 542291761582 & 103 & 652021378808 & 143 & 676166974426\\
24 & 255779464606 & 64 & 1118566817183 & 104 & 347023836636 & 144 & 736507281912\\
25 & 720108874443 & 65 & 542326334671 & 105 & 1377310608996 & 145 & 706530314192\\
26 & 617516839615 & 66 & 1278162950931 & 106 & 610461715104 & 146 & 581165299308\\
27 & 465700248000 & 67 & 618114265525 & 107 & 874154628868 & 147 & 5276071714148\\
28 & 715585432812 & 68 & 766809981348 & 108 & 264142448144 & 148 & 2643550739558\\
29 & 1482967536426 & 69 & 1369524128609 & 109 & 1369524128609 & 149 & 3607718249282\\
30 & 451696258296 & 70 & 728617527999 & 110 & 372699692218 & 150 & 2227364770841\\
31 & 728857315908 & 71 & 726799330926 & 111 & 721274949941 & 151 & 1125014606729\\
32 & 704735395618 & 72 & 1210753867510 & 112 & 1183740987404 & 152 & 2025616102588\\
33 & 347212673106 & 73 & 955358130060 & 113 & 388455096872 & 153 & 1776255312031\\
34 & 921254501510 & 74 & 209410626542 & 114 & 2650373378946 & 154 & 596281032778\\
35 & 913255188514 & 75 & 967895418038 & 115 & 1015354182030 & 155 & 419583263768\\
36 & 306628397688 & 76 & 1497895417012 & 116 & 716928899496 & 156 & 1125014606729\\
37 & 712073900188 & 77 & 482819399640 & 117 & 1017286657609 & 157 & 332343732936\\
38 & 1198423260540 & 78 & 1202370080434 & 118 & 427599690212 & 158 & 2025616102588\\
39 & 612076728015 & 79 & 632759531960 & 119 & 687777826819 & 159 & 974118613240\\
40 & 884538517819 & 80 & 792312038485 & 120 & 1295945493968 & Total & 149364113290700\\
\hline
\end{tabular}
\end{table}

\begin{sidewaystable}
\centering
\caption{The permanent values of $123$ sub-matrices for $C_{60}$. \label{T14}}
\begin{tabular}{cccccc}
\hline
$No.$ & $per$ & $No.$ & $per$ & $No.$ & $per$\\
\hline
1 & -108132060044208 - 17236666752811.5i & 43 & 426568556602.3 + 241672686883.4i & 85 & 540889261725.6 - 104622128811.5i\\
2 & 32439176143560.7 + 1834294235804.6i & 44 & 107807473965.4 + 73810357204.4i & 86 & -817143586156.8 + 47166381251.1i\\
3 & 28651718827889.8 + 1357573987458.8i & 45 & {\tiny -269465247652081.6e+14 + 706095638266417.0e+13i} & 87 & -994806975892.5 + 93886331421.2i\\
4 & 6789299586006.0 - 1688236535.6i & 46 & {\tiny -661668903812826.1e+27 + 492728489717893.1e+27i} & 88 & 8109989785555.2 - 2045825961592.6i\\
5 & -2173724934049.1 - 234478873974.3i & 47 & {\tiny -570045180517617.4e+28 + 593490587519625.4e+28i} & 89 & 1388436024137.2 - 254706472245.7i\\
6 & {\tiny 205553118230115.7e+23 - 141243991818827.9e+22i} & 48 & 2005516079797.5 - 264097579607.2i & 90 & 946931667201.5 - 41619802367.6i\\
7 & -4208339700022.5 + 48197250732.1i & 49 & -6923809018625.4 - 925864779886.0i & 91 & 140310887337.4 - 101777728465.5i\\
8 & -3715786162910.6 - 426114665504.4i & 50 & -1467766896102.5 - 295130882950.8i & 92 & 2015007303507.5 + 37289275154.3i\\
9 & 5062860050.5 + 1372978026.9i & 51 & -1385153550914.2 - 114635561390.5i & 93 & 676228873833.0 - 55610438973.4i\\
10 & 5422484098412.9 - 942316127646.7i & 52 & -947000498530.3 - 35797287385.8i & 94 & 37837324337277.8 - 1269322916070.9i\\
11 & 13196446453306.6 - 652755644317.2i & 53 & -136871010093.6 - 11543917443.7i & 95 & 14744473937369.5 - 1210958323675.0i\\
12 & 4706167841951.5 + 854731960557.2i & 54 & 5815633975.2 - 3473208780.2i & 96 & 5729268103836.0 + 42724101089.7i\\
13 & 2046379439837.8 - 59442924832.4i & 55 & -116800600530.1 + 29537465350.0i & 97 & 702701469931.9 + 48964516611.8i\\
14 & 1410002078478.6 + 95635056612.4i & 56 & -2208745388113.2 + 335084956770.2i & 98 & 1155277939675.6 + 148629934592.4i\\
15 & 1617597373807.2 - 270511638657.2i & 57 & 535498735978.0 - 79221341491.8i & 99 & 356756381159.4 + 19117604919.0i\\
16 & 5670038850189.0 - 1101990793674.7i & 58 & -46070027008005.2 + 5929984393672.9i & 100 & 6549241518036.4 - 59483813012.0i\\
17 & 1461009383950.5 - 569793948123.5i & 59 & -11243051367865.9 + 520387307092.4i & 101 & 928877327413.1 + 112834748074.2i\\
18 & -1929245154069.7 + 413889668728.7i & 60 & 5198433030736.2 - 297381728075.6i & 102 & 753424498572.8 - 212580395371.2i\\
19 & -12855617174845.0 + 253069807489.3i & 61 & 1387004773127.8 + 58135044793.9i & 103 & 170765066948.4 + 35221340229.0i\\
20 & -3963324935926.1 - 349976755162.1i & 62 & 988155323940.3 - 145293280529.6i & 104 & 1148546296973.9 - 96002491933.5i\\
21 & -16043983372582.8 - 738837585232.1i & 63 & 1987844225942.7 + 275773954315.2i & 105 & 318327255660.9 - 98255335916.4i\\
22 & 1616424380561.5 - 348664794337.4i & 64 & -5137551856506.9 + 212228141173.9i & 106 & 417200651281.2 - 26347380961.5i\\
23 & -569574088601.2 - 35581542923.8i & 65 & -867124998758.5 - 36317187545.5i & 107 & 5024288099014.3 - 18055323710.0i\\
24 & 208761634761.6 + 100224643143.9i & 66 & -881718222829.8 - 55460487636.4i & 108 & -1551343571116.4 + 63433450818.7i\\
25 & -2200105202590.9 - 193752163463.9i & 67 & -2371097760173.1 - 303687098188.3i & 109 & 669483634510.3 + 51886287177.4i\\
26 & -700444445654.6 + 17906258599.9i & 68 & -671850169543.6 - 215518096819.5i & 110 & 552775115082.4 - 127131168691.7i\\
27 & -1457262659340.6 + 97692598693.7i & 69 & 3295485506358.3 + 316379606154.8i & 111 & 165662792868.9 - 22937653807.9i\\
28 & -4272295160158.9 - 818731359827.5i & 70 & -4328333841250.1 + 272826336708.3i & 112 & 200060545615.8 + 60234117810.4i\\
29 & -1727536315959.6 + 99086669151.5i & 71 & 1878975036662.2 - 482764646653.9i & 113 & 1257142874363.4 - 190010410554.8i\\
30 & 435168291995.2 - 233828091.5i & 72 & 348166392620.7 - 19625110968.7i & 114 & -232594984598.9 + 51677492433.3i\\
31 & -200943978554.4 + 98296317364.4i & 73 & 719383143994.6 + 110029768817.6i & 115 & 302577563322.0 - 73839010357.0i\\
32 & 6135826550.3 - 5964846837.2i & 74 & 1739679686167.9 + 52621999591.9i & 116 & 372434227966.9 - 39687783432.4i\\
33 & 41251773751.7 + 12060793948.7i & 75 & 4304376060302.4 - 811532609944.6i & 117 & 12882964806795.9 + 2212041426209.0i\\
34 & -10073588432.7 - 917393934.4i & 76 & 1553023289221.0 + 246572995110.8i & 118 & 1674327862222.5 + 237454637881.9i\\
35 & 2308788385165.9 + 5083792535656.3i & 77 & -150163473464.0 - 4383160522.4i & 119 & 2413900722675.0 - 112334487535.5i\\
36 & 8834947401589.9 + 572758118440.2i & 78 & -14892023440279.6 + 1690606913584.6i & 120 & 406823127969.1 + 43187852505.6i\\
37 & 5398524831188.1 + 1012206133038.6i & 79 & -3268390974156.7 + 663085472835.0i & 121 & 2651565412357.1 + 511765637640.9i\\
38 & 2562921080124.5 + 294888019848.2i & 80 & -1501750651439.6 - 53619370473.0i & 122 & 337999903586.1 + 38393644710.3i\\
39 & 916011724505.2 - 37527124281.4i & 81 & -537787336683.5 - 82992336064.9i & 123 & 1003141670015.9 + 137179900724.3i\\
40 & {\tiny -718428632724154.8e+13 + 155895693830716.6e+13i} & 82 & -552003362140.4 - 99045880290.3i & Total & {\tiny -729320933310243.2e+28 + 631208207405806.6e+28i}\\
41 & {\tiny -931109179425221.2e+27 - 115550878416171.2e+27i} & 83 & -3548113350360.1 + 640548481115.9i\\
42 & 1209454910810.9 - 195013751259.9i & 84 & 498121576738.7 - 48173604015.8i\\
\hline
\end{tabular}
\end{sidewaystable}

\end{document}